\documentclass[a4paper,12pt]{amsart}

\usepackage[english]{babel}
\usepackage{amsthm}
\usepackage{amssymb}
\usepackage{hyperref}

\newcommand{\Set}[1]{\left\{\, #1 \,\right\}}

\newcommand{\Span}[1]{\langle\, #1 \,\rangle}
\newcommand{\BigSpan}[1]{\left\langle\, #1 \,\right\rangle}

\newcommand{\Size}[1]{\lvert #1 \rvert}

\DeclareMathOperator{\im}{im}

\DeclareMathOperator{\GL}{GL}
\DeclareMathOperator{\GF}{GF}
\DeclareMathOperator{\Aut}{Aut}
\DeclareMathOperator{\End}{End}

\DeclareMathOperator{\Hom}{Hom}

\renewcommand{\phi}[0]{\varphi}
\renewcommand{\theta}[0]{\vartheta}
\renewcommand{\epsilon}[0]{\varepsilon}

\numberwithin{dummy}{section}
\numberwithin{equation}{section}

\theoremstyle{definition}

\begin{document}

\date{12 September 2013, 13:36 CEST --- Version 5.04}

\title[Abelian automorphism groups]
{A module-theoretic approach\\
to abelian automorphism groups}

\author{A.~Caranti}

\address[A.~Caranti]%
 {Dipartimento di Matematica\\
  Universit\`a degli Studi di Trento\\
  via Sommarive 14\\
  I-38123 Trento\\
  Italy} 

\email{andrea.caranti@unitn.it} 

\urladdr{http://science.unitn.it/$\sim$caranti/}

\begin{abstract}
  There are  several examples in the literature  of finite non-abelian
  $p$-groups whose automorphism group  is abelian.  For some time only
  examples  that were special  $p$-groups were  known, until  Jain and
  Yadav~\cite{JY12}  and Jain,  Rai  and Yadav~\cite{JRY}  constructed
  several non-special  examples. In  this paper we  show how  a simple
  module-theoretic approach allows the construction of non-special
  examples, starting from special ones constructed by several authors,
  while at the same time avoiding further direct calculations.
\end{abstract}

\keywords{finite $p$-groups, automorphisms, central automorphisms}

\thanks{The author gratefully acknowledges the support of the
  Department of Mathematics, University of Trento}
  
\subjclass[2010]{20D15 20D45}

\maketitle

\thispagestyle{empty}

\section{Introduction}

There  are several examples  in the  literature of  finite, non-abelian
groups with an  abelian automorphism group, starting with  the work of
G.A.~Miller  \cite{Miller13, Miller14}  of a  century ago.  Clearly in
such a group $G$ the central  quotient $G / Z(G)$, which is isomorphic
to the group of inner automorphisms of $G$, is abelian, so that $G$ is
of nilpotence class $2$.

Until recently,  all known examples of  finite, non-abelian $p$-groups,
for  $p$  an  odd  prime,  with an  abelian  automorphism  group  were
\emph{special} $p$-groups, satisfying
\begin{equation}\label{eq:special_special}
  \begin{aligned}
    &G^{p} \le G' = \Phi(G) = Z(G) = \Omega_{1}(G),
        \\&\Size{G / G'} = p^{n} = \Size{G^{p}}, 
        \quad
        \Size{G'} = p^{\binom{n}{2}},
  \end{aligned}
\end{equation}
for some $n \ge 4$. (See Section~\ref{sec:modules} for a generalization.)

In such a group $G$, the group $\Aut_{c}(G)$ of central automorphisms,
that is, those  automorphisms that induce the identity  on the central
quotient $G  / Z(G)$,  is an elementary  abelian group of  order $p^{n
  \binom{n}{2}}$.     Several    authors    have    constructed
examples~\cite{HeLie,  JoKo, Ea, H79,  CaLe, Ca83,  Ca85, Mo94,  Mo95} of
groups  $G$ as in~\eqref{eq:special_special},  in which  the relations
are suitably chosen as to  force all automorphisms to be central, that
is, $\Aut(G) = \Aut_{c}(G)$. It follows that $\Aut(G)$ is (elementary)
abelian.

In Section~\ref{sec:modules} of this paper we show that this situation
allows for a simple module-theoretic translation. For previous usages
of this technique, see~\cite{HeLie, DaHe, Ca83, GPS}. Namely, both $V =
G/G'$ and 
$G'$ can be regarded as vector spaces over the field $\GF(p)$ with $p$
elements.   The action of  $\GL(V)$ on  $V =  G/G'$ induces  a natural
action on $G'$. Moreover, the map  $a \wedge b \mapsto [a, b]$ extends
to an  isomorphism of vector spaces  from   the exterior square
$\Lambda^{2}   V$ to $G'$,  and   this   map  is   also   an
isomorphism   of 
$\GL(V)$-modules.  Consider the  linear map $f : V  \to \Lambda^{2} V$
given  by  $(a  G')f =  a^{p}$,  where  we  are identifying  $G'$  and
$\Lambda^{2} V$.  (And we write maps on the right.) Then the statement
that  all  automorphisms of  $G$  are  central  is equivalent  to  the
statement that
\begin{equation*}
  \Set{ \alpha \in \GL(V) : \alpha \circ f = f \circ \hat{\alpha} }
  =
  \Set{1},
\end{equation*}
where $\hat{\alpha}$ is the automorphism of $\Lambda^{2} V$ induced by
$\alpha$. 

Mahalanobis had conjectured~\cite{Ma08} that all finite $p$-groups
whose automorphism group is abelian are special. However, Jain and
Yadav~\cite{JY12} and Jain, Rai and Yadav~\cite{JRY} have constructed
several examples to show that this is not the case. These examples
require ingenious and extensive calculations to determine the
automorphism groups.

In  the  rest   of  this  paper  we  show   how  the  module-theoretic
reformulation  of  Section~\ref{sec:modules}  allow  us  to  construct
examples  similar to  those  of~\cite{JY12, JRY}  within a  conceptual
framework, which does not  require further calculations.  In fact, we
start  from  known  examples of  the  form~\eqref{eq:special_special},
modify them  in simple ways,  and use the module-theoretic  setting to
show that  the required  properties of the  automorphism group  of the
modified  groups   follow  rather  directly   from  the  corresponding
properties  of  the  original   groups,  without  resorting  to  extra
calculations.   This  allows  us   to  construct  examples  of  finite
$p$-groups $G$ with the following properties:
\begin{enumerate}
\item $\Aut(G) = \Aut_{c}(G)$ is non-abelian, and $G$ does not have
  abelian, nontrivial direct factors (Section~\ref{sec:Zurek});
\item $\Aut(G) = \Aut_{c}(G)$ is elementary abelian, and $G' = \Phi(G)
  < Z(G)$, so that $G$ is not special (Section~\ref{sec:bigger_centre});
\item $\Aut(G) = \Aut_{c}(G)$ is abelian, and $G' < \Phi(G) = Z(G)$,
  so that   $G$ is not special (Section~\ref{sec:bigger_frattini}).
\end{enumerate}

Recall  that   a  non-abelian  group   is  said  to   be  \emph{purely
  non-abelian}  if it  has only  trivial abelian  direct  factors. Non
purely   non-abelian  examples   $G$  with   non-abelian   $\Aut(G)  =
\Aut_{c}(G)$    were    constructed    by   Curran~\cite{Cur82}    and
Malone~\cite{Mal84}.   Glasby~\cite{Gla86}   constructed  an  infinite
family of  purely non-abelian, finite $2$-groups  $G$ with non-abelian
$\Aut(G)  =  \Aut_{c}(G)$.   Jain  and  Yadav~\cite[Theorem~B]{JY12}
have constructed 
examples as  in~(1), that is,  purely non-abelian $p$-groups  $G$ with
non-abelian $\Aut(G) = \Aut_{c}(G)$, for $p$ an odd prime.

Examples as in~(2) have been constructed in~\cite[Theorem~B]{JRY}. 

Examples as in~(3)
have been constructed in~\cite[Theorem~A, Lemma~2.1]{JY12}
and~\cite[Theorem~C]{JRY}, the latter groups having elementary
abelian automorphism group. 

We  are  very grateful  to  Carlo  Scoppola,  Manoj K.~Yadav,  Yassine
Guerboussa and  Sandro Mattarei for their useful  observations. We are
indebted   to  the   referee  for   several  important   comments  and
suggestions.

\section{Notation}

Throughout the paper, $p$ denotes an odd prime.

Our iterated commutators are left-normed, $[a, b, c] = [[a,b],c]$.

A (non-abelian) $p$-group $G$ is said to be \emph{special} if $G^{p}
\le G' = Z(G)$. This implies that $\Phi(G) = G^{p} G' = Z(G)$, and that
$G'$ is elementary abelian. 

A  non-abelian  group   is  said  to   be  \emph{purely
non-abelian}  if it  has only  trivial abelian  direct  factors.

We will be writing maps on the right, and group morphisms as
exponents, so  if $g \in G$ and 
$\phi \in \Aut(G)$, then $g^{\phi}$ denotes the action of $\phi$ on
$g$. We also write $g^{\phi + \psi} = g^{\phi} g^{\psi}$, for $\phi,
\psi \in \Aut(G)$, and in other similar situations.


\section{The module-theoretic approach}
\label{sec:modules}

Suppose we have a group $G$ as in~\eqref{eq:special_special}, for an
odd prime $p$, and some $n
\ge 4$. Then $G$ admits a presentation of the form
\begin{equation}\label{eq:presentation}
  \begin{aligned}
    G
    =
    \Span{x_{1}, \dots, x_{n}
      :\
      &\text{$[x_{i}, x_{j}, x_{k}] = 1$ for all $i, j, k$, }
      \\&\text{$[x_{i}, x_{j}]^{p} = 1$ for all $i, j$, }
      \\&\text{$x_{i}^{p} = \prod_{j < k} [x_{j}, x_{k}]^{c_{i,j,k}}$ for all $i$}
    },
  \end{aligned}
\end{equation}
where the $c_{i,j,k}$ are such that the $x_{i}^{p}$ are independent in
$G'$,  for $i  = 1,  \dots,  n$. Note  that we  have explicitly  added
relations that state that $G'$ is of exponent $p$; these are redundant
here, as  they can be deduced  from $1 = [x_{i}^{p},  x_{j}] = [x_{i},
  x_{j}]^{p}$,  but will  play a  role in  the generalizations  in the
later sections.  Note also that the relations imply that $G^{p} \le G'
= Z(G)$, so that $G' =  \Phi(G)$, and also, as we just said, $(G')^{p}
= 1$.  Moreover $\Size{G/G'} =  \Size{G^{p}} = p^{n}$ and $\Size{G'} =
p^{\binom{n}{2}}$, as indeed required in~\eqref{eq:special_special}.


Consider the vector space  $V = G/G'$ over $\GF(p)$.
Let $\phi \in \Aut(G)$, and let $\alpha$ be the automorphism of
$V = G / G'$ induced by $\phi$, so that we may regard $\alpha$ as an
element of $\GL(V)$. The action of $\phi$ on $G'$ is completely
determined by $\alpha$. In fact, since $G'\le Z(G)$ we have
\begin{equation*}
  [a, b]^{\phi} = [a^{\phi}, b^{\phi}] = [(a G')^{\alpha}, (b G')^{\alpha}].
\end{equation*}
We denote by
$\hat{\alpha}$ the automorphism of $G'$ induced by $\alpha$ (or
$\phi$). 
The
map $a \wedge  b \mapsto [a,  b]$ extends to an  isomorphism of vector
spaces from  the  exterior square $\Lambda^{2} V$ to $G'$, and this map
is also an isomorphism of  $\GL(V)$-modules.  Thus we can associate to
such a group the linear map defined by
\begin{align*}
  f :\ &V \to \Lambda^{2} V\\
       &x_{i} \mapsto \sum_{j, k} c_{i, j, k} \cdot x_{j} \wedge x_{k}.
\end{align*}
 
In the papers~\cite{HeLie,  JoKo, Ea, H79, CaLe, Ca83,  Ca85, Mo94, Mo95}
examples are constructed of groups of the form~\eqref{eq:presentation}
in which all automorphisms are central (see below for variations).


In all of these examples, to prove that all automorphisms are central one
uses the fact that a map $\Set{x_{1}, \dots, x_{n} } \to G$ extends to
an endomorphism  of $G$ if and  only if it preserves  the $p$-th power
relations  in~\eqref{eq:presentation}.  

This is clear in one direction. If $\phi \in \End(G)$, and $y_{i} =
x_{i}^{\phi}$, then the $y_{i}$ satisfy
\begin{equation}\label{eq:rels-for-y}
  y_{i}^{p} = \prod_{j < k} [y_{j}, y_{k}]^{c_{i,j,k}}, \quad
  \text{for $i = 1, \dots, n$.}
\end{equation}
Conversely, suppose we  have a map $\Set{x_{1}, \dots,  x_{n} } \to G$
such  that  $x_{i}  \mapsto   y_{i}$,  and  assume  that  the  $y_{i}$
satisfy~\eqref{eq:rels-for-y}.  Consider  the free  group  $F$ in  the
variables $z_{1},  \dots, z_{n}$  in the variety  defined by  the laws
$[a, b, c], [a, b]^{p}$, that  is, the variety of groups of nilpotence
class  at  most  $2$,  with  derived  subgroup  of  exponent  dividing
$p$. Clearly $G$ belongs to this  variety, and the epimorphism $F \to
G$ determined by $z_{i}  \mapsto x_{i}$ has  kernel $N$ 
generated by the
\begin{equation*}
  z_{i}^{-p} \prod_{j < k} [z_{j}, z_{k}]^{c_{i,j,k}}, \quad \text{for
    $i = 1, \dots , n$.}
\end{equation*}
Identify $G$ with $F/N$, by taking $x_{i} = z_{i} N$. Since the
$y_{i}$ satisfy~\eqref{eq:rels-for-y}, the  morphism $F
\to G$ which sends $z_{i} \mapsto y_{i}$ has $N$ in its kernel, and
thus induces a morphism $G = F/N \to G$ which sends $x_{i} = z_{i} N
\mapsto y_{i}$. This is the required endomorphism of $G$.

Thus one  starts with a  map $\Set{x_{1}, \dots,  x_{n} } \to  G$, and
assuming it extends  to $\phi \in \Aut(G)$, one  computes the value of
$x^{p}  \phi$,  as $x$  ranges  over  a basis  of  $V  =  G/G'$, in  two
ways. First,  since $p > 2$,  the $p$-th power map  induces a morphism
$G/G' \to G'$.  Thus, in the notation just introduced, $(x^{p})^{\phi}
= (x^{\phi})^{p}  = ((x G')^{\alpha})^{p}$, and this  yields one value
of  $(x^{p})^{\phi}$.   Another  value  comes from  $(x^{p})^{\phi}  =
(x^{p})^{\hat{\alpha}}$,  as  $x^{p}  \in  G'$.   Equating  these  two
values, one  obtains (quadratic) equations for the  entries of $\alpha
\in  \GL(V)$ as a  matrix with  respect to  a suitable  basis of  $V =
G/G'$, and  these equations are exploited  to show that  $\alpha = 1$,
that is $\phi \in \Aut_{c}(G)$.

As we said in the Introduction, in module-theoretic terms, 
the fact that $\Aut(G) = \Aut_{c}(G)$
is thus seen to be equivalent to saying that the only element
of $\GL(V)$ that commutes with $f$ is the identity,
\begin{equation*}
  \Set{ \alpha \in \GL(V) : \alpha \circ f = f \circ \hat{\alpha} }
  =
  \Set{1}.
\end{equation*}

It might be remarked here that in~\cite{HeLie, DaHe, Ca85, GPS} a more
general fact has been used. This states
that there is 
an automorphism $\phi$ of $G$ which induces the automorphism $\alpha$
of $V = G/G'$ if and only if 
\begin{equation}\label{eq:the-switch}
  \alpha \circ f = f \circ
\hat{\alpha}.
\end{equation}
The direct implication is clear. Conversely, if~\eqref{eq:the-switch}
holds, choose 
elements $y_{i} \in G$ such that 
$(x_{i} G') \alpha = y_{i} G'$. Because of~\eqref{eq:the-switch}, the
$y_{i}$ satisfy~\eqref{eq:rels-for-y}, and thus the map $\Set{x_{1},
  \dots, x_{n} } \to G$ that sends $x_{i} \mapsto y_{i}$ extends to an
automorphism $\phi$ of $G$, which induces $\alpha$.

The above setup admits a slight generalization
. In some examples
of groups with $\Aut(G) = \Aut_{c}(G)$,
the presentation~\eqref{eq:presentation} includes some relations within
$G'$, so that some products of commutators vanish. However, we require
that the condition $Z(G) = \Phi(G)$ still holds,
that is, no generator is made central by adding these commutator
relations, and also that the relation
$\Size{G/G'} = \Size{G^{p}}$ is preserved. Some examples, for
instance in~\cite{Mo94}, make do without the last requirement, but we
will need it for the examples in the following sections.

So we have presentations of the form
\begin{equation}\label{eq:presentation_K}
  \begin{aligned}
    G
    =
    \Span{x_{1}, \dots, x_{n}
      :\
      &\text{$[x_{i}, x_{j}, x_{k}] = 1$ for all $i, j, k$, }
      \\&\text{$[x_{i}, x_{j}]^{p} = 1$ for all $i, j$, }
      \\&\text{$\prod_{j < k} [x_{j}, x_{k}]^{d_{l,j,k}} = 1$
        for $l = 1, \dots , t$, }
      \\&\text{$x_{i}^{p} = \prod_{j < k} [x_{j}, x_{k}]^{c_{i,j,k}}$
        for all $i$}
    },
  \end{aligned}
\end{equation}
for some $t$, with $p^{n} = \Size{G/G'} = \Size{G^{p}}$, and $Z(G) =
\Phi(G)$.   

When determining the automorphisms of $G$, starting with $\alpha \in
\GL(V)$ (where $V = G/G'$), we have to make sure that the
automorphism $\hat{\alpha}$ leaves invariant the subspace
\begin{equation}\label{eq:K}
  K =
  \BigSpan{\sum_{j < k} d_{l,j,k} \cdot x_{j} \wedge x_{k} : l = 1, \dots, t}
\end{equation}
of
$\Lambda^{2} V$, which corresponds to the extra relations within $G'$.

Thus in  the module-theoretic  setting we have 
\begin{enumerate}
\item 
a vector space  $V$ of
dimension at least  $4$ over $\GF(p)$,
\item a subspace  $K$ of $\Lambda^{2}
V$ such that $v \wedge V \not\subseteq
K$ for $0 \ne v \in V$  (this is a translation of the requirement
$Z(G) = \Phi(G)$), and 
\item an \emph{injective}  linear map 
\begin{equation}\label{eq:f}
  \begin{aligned}
  f :\ &V  \to \Lambda^{2}  V / K
     \\&v_{i} \mapsto \sum_{j, k=1}^{n} c_{i,j,k} \cdot v_{j} \wedge v_{k},
  \end{aligned}
\end{equation}
where $v_{1},  \dots, v_{n}$ is  a basis of  $V$.  
\end{enumerate}
We say that  such a
triple  $(V,  K,  f)$  is  a \emph{Trivial  Automorphism  Triple},  or
\emph{TAT}, if in addition to (1)-(3) above, we have
\begin{enumerate}
\item[(4)]
$
  \Set{ \alpha \in \GL(V) : K^{\hat{\alpha}} = K,
    \alpha \circ f = f \circ \bar{\alpha} }
  =
  \Set{1}.
$
\end{enumerate}
Here $\hat{\alpha}$ and $\bar{\alpha}$ denote the automorphisms
induced by $\alpha$ respectively on $\Lambda^{2}  V$
and, since $K$ is left invariant by $\hat{\alpha}$, on
$\Lambda^{2}  V / K$.
Thus every group $G$ of the
form~\eqref{eq:presentation_K} for which  $\Aut(G) = \Aut_{c}(G)$
yields a TAT. 


Conversely, given a TAT 
$(V, K, f)$, we write succinctly 
\begin{equation}\label{eq:presentation_compact}
  G = \Span{x : x^{p} = x f, K}
\end{equation}
for  the  presentation~\eqref{eq:presentation_K},   where  $K$  is  as
in~\eqref{eq:K} and $f$ as in~\eqref{eq:f}. Our discussion shows that
for such a group $G$ we 
have $\Aut(G) = \Aut_{c}(G)$. 


\section{A group $G$ with all automorphisms central, but $\Aut(G)$  non-abelian}
\label{sec:Zurek}

In this  section we apply the  techniques of Section~\ref{sec:modules}
to construct an  example of a finite, purely  non-abelian $p$-group $G$
such that $\Aut(G) = \Aut_{c}(G)$ is non-abelian. 

Let  us first  review  why $\Aut(G)  =  \Aut_{c}(G)$ is  abelian in  a
group~\eqref{eq:presentation_compact}, where $(V, K, f)$ is a
TAT. (This would also follow from~\cite[Theorem 4]{AdYen}, reported
as~\cite[ Theorem~2.7]{JRY}.) 

An observation of Guerboussa and
Daoud~\cite{GuDa}, which is a specialization of results of
H.~Laue~\cite{Laue85}, comes handy. This states that if $G$ is purely
non-abelian (a condition that is satisfied here, 
as we require $Z(G) \le \Phi(G)$), then $\Aut_{c}(G)$ under
composition is isomorphic to 
$\Hom(G, Z(G))$ under the circle operation $\gamma \circ \delta =
\gamma + \delta + \gamma \delta$, where $\gamma \delta$ denotes
composition. (Recall that $g^{\gamma + \delta} = 
g^{\gamma} g^{\delta}$, and note that $\gamma + \delta = \delta +
  \gamma$, as $g^{\gamma}, g^{\delta} \in Z(G)$.)
This is readily seen, as if $\phi$ is a central automorphism, then the
map $\gamma$ given by
 $x^{\gamma} = x^{-1} x^{\phi} = [x, \phi]$ is in $\Hom(G, Z(G))$.

In the case
of~\eqref{eq:presentation_compact}, we have $G' = Z(G)$, so if
$\gamma, \delta \in \Hom(G, Z(G))$, then $\im(\gamma), \im(\delta) \le
Z(G) = G' \le
\ker(\gamma), \ker(\delta)$, so that $\gamma \delta = 0 = \delta \gamma$,
and $\gamma \circ \delta = \gamma + 
\delta + \gamma \delta = \gamma + \delta = \delta + \gamma = \delta
\circ \gamma$. This 
shows that $\Aut(G)$ is abelian.

We now exploit an idea of Zurek~\cite{Zu}. We start with a TAT $(V, K, f)$
and construct a group as 
the following variation on~\eqref{eq:presentation_compact}
\begin{equation*}
  G = \Span{x : x^{p^{2}} = x f, K}.
\end{equation*}
By this we mean a presentation like~\eqref{eq:presentation_K}, where
we have replaced the $p$-th powers in the last line of relations with
$p^{2}$-th powers. Note that $G$ is still in the variety of groups of
class $2$, and derived
subgroup of exponent $p$. Also
$G' < \Phi(G) = Z(G)$, and the $p^{2}$-th power map is a morphism
$G/\Phi(G) \to G'$ which, with the proper identifications,  is the
same as $f : V \to \Lambda^{2} / K$. The same arguments of the
previous section apply 
to show that the fact that $(V, K, f)$ is a TAT implies that all
automorphisms of $G$ are 
central. 

To show that $\Aut(G)$ is non-abelian, we use again the observation of
Guerboussa and Daoud: let $n = \dim(V)$, and let $x_{1}, \dots, x_{n}$
be a basis for $V$. Define $\gamma, \delta \in \Hom(G, Z(G))$ by
\begin{equation*}
  \begin{aligned}
      &x_{1}^{\gamma} = x_{1}^{p}, \quad x_{i}^{\gamma} = 1\text{, for $i > 1$}
      \\&x_{1}^{\delta} = x_{2}^{p}, \quad x_{i}^{\delta} = 1\text{, for $i > 1$}.
  \end{aligned}
\end{equation*}
Then $x_{1}^{\gamma \delta} = (x_{1}^{p})^{\delta} = x_{2}^{p^{2}} \ne 1 =
 (x_{2}^{p})^{\gamma} = x_{1}^{\delta \gamma}$, so that $\Aut(G) =
\Aut_{c}(G)$ is non-abelian.

\section{Groups with $\Aut(G) = \Aut_{c}(G)$ abelian, and $G' = \Phi(G) < Z(G)$}
\label{sec:bigger_centre}

In this section we construct an example of a finite $p$-group $G$ with
$\Aut(G) = \Aut_{c}(G)$ elementary abelian, and $G' = \Phi(G) <
Z(G)$, so that $G$ is not special. 

We start with a TAT $(V, K, f)$ such that $K + V f$ is a proper subspace of $\Lambda^2 V$ 
(for instance $K = \{ 0 \}$ will do), and construct first a group 
\begin{equation}\label{eq:mat1}
  H = \Span{x : x^{p} = x f, K}
\end{equation}
as in~\eqref{eq:presentation_K}. 

Then we construct a group $G$ as a central product $H \Span{z}$, where $z$ has
order $p^{2}$, and we amalgamate $z^{p} \in H' \setminus H^{p}$. (The extra condition on the TAT in fact implies that $H^p < H'$.)

We have $W = G/G' \cong V \oplus \Span{w}$, where $V = H/H'$ and $w$
is the image of $z$. 
Consider $\phi \in \Aut(G)$. Clearly $\Span{w} = Z(G) \, G' / G'$ is left
invariant by $\phi$. Write the action of $\phi$ on $W$ in matrix form
\begin{equation}
  \begin{bmatrix}
    \alpha & \lambda\\
    0      & \mu
  \end{bmatrix},
\end{equation}
with respect  to a basis of  $V$, extended with $w$.  Here $\alpha \in
\GL(V)$, $\mu \in \GF(p)$, and  $\lambda$ is a column vector of length
$n$.  (Since our  automorphisms act on the right,  elements of $V$ and
$W$ are row vectors.)

Consider  the  embedding  $\iota  :  V \to  V  \oplus  \Span{w} = W$.  The
composition
\begin{equation*}
    V \overset{\iota}{\to}  V \oplus \Span{w} \overset{\cdot^{p}}{\to} 
    \Lambda^{2} V
\end{equation*}
equals $f$. Moreover, the action of $\phi$ on $G' \cong
\Lambda^{2} V$ is completely determined by $\alpha$, as $z \in
Z(G)$. Therefore we have
\begin{equation*}
  v^{\alpha} f = (v f)^{\hat{\alpha}},
\end{equation*}
so that $\alpha = 1$, as $(V, K, f)$ is a TAT. This implies that
$\hat{\alpha} = 1$ on $\Lambda^{2} V$, and since the $p$-th power map
$V \oplus \Span{w} \to \Lambda^{2} V$ is injective, also $\lambda = 0$
and $\mu = 1$. It follows that the matrix in~\eqref{eq:mat1} is the
identity, and thus $\phi$ induces the identity modulo $G'$. This shows
that $\Aut(G) = \Aut_{c}(G)$.

We have shown that all automorphisms of $G$ are trivial modulo $G' <
Z(G)$. Since $G'$ is elementary abelian, $\Aut(G)$ is also elementary
abelian, as it is isomorphic to $(\Hom(G/G', G'), \circ)$.

\section{Groups with $\Aut(G) = \Aut_{c}(G)$ abelian, and $G'  < \Phi(G) = Z(G)$}
\label{sec:bigger_frattini}

In this section we construct a finite $p$-group $G$ with
$\Aut(G) = \Aut_{c}(G)$ abelian and $G'  < \Phi(G) = Z(G)$; in particular
$G$ is 
not special. 

Let  us  start  from  a  TAT  $(V,  \Set{0},  f)$,  with  $\Size{V}  =
p^{n}$. Define a special $p$-group
\begin{equation*}
  H = \Span{x : x^p = x f}
\end{equation*}
as in~\eqref{eq:presentation_compact}, with $K = \Set{0}$.


As seen in  Section~\ref{sec:Zurek}, if $V = H/H'$,  then the group of
central  automorphisms of $H$  is isomorphic  to $\Hom(H,  Z(H)) \cong
\Hom(V,  \Lambda^{2}  V)$  under  the  circle  operation.  The  latter
contains the subspace
\begin{equation*}
\Set{x  \mapsto x \wedge v : v \in V}
\end{equation*}
of the \emph{inner} maps,  which correspond to the inner automorphisms
of $H$. Since $v \wedge v = 0$  for all $v \in V$, none of these inner
maps is injective. Therefore  if we choose an \emph{injective} $\gamma
\in \Hom(H, Z(H))$ (for instance, $x \mapsto x^{p}$ will do), this
will not be inner. Let $m > 1$, and consider the 
extension  $G$  of $H$  by an  element $y$  such that  $[x, y]  = x
\gamma$, for  $x \in H$, and  $y^{p^{m}} \in H'  \setminus H^{p}$. The
choice of $\gamma$ ensures that $Z(G) \le \Phi(G)$, and actually $Z(G)
= \Phi(G)$.

We have $G/\Phi(G) \cong \Span{u} \oplus V$, where $u$ is the image of
$y$.   Let   $\phi  \in   \Aut(H)$.   Since  $V$   is  the   image  
in $\Span{u} \oplus
V$ of
$\Omega_{2}(G)$, the  automorphism $\phi$ induces  on $\Span{u} \oplus
V$ a linear map which can be described by the matrix
\begin{equation}\label{eq:mat2}
  \begin{bmatrix}
    \tau & \sigma\\
    0      & \alpha
  \end{bmatrix},
\end{equation}
with respect  to a basis which  begins with $u$, and  continues with a
basis of $V$. Here  $\alpha \in \GL(V)$, $\tau \in \GF(p)$, and $\sigma$ is a row vector of length $n$. Since $G' = H'$, we have
that $\alpha$ alone determines  the action $\hat{\alpha}$ of $\phi$ on
$G'$.

As in the previous section, the composition
of the embedding $V \overset{\iota}{\to}  \Span{u} \oplus V$ with the $p$-th power map
equals $f$.  Since $(V,  \Set{0}, f)$  is a TAT,  we obtain  $\alpha =
\hat{\alpha} = 1$.  In particular $\phi$ acts trivially  on $G' \cong
\Lambda^{2} V$. This 
also implies $\tau = 1$.

Let $x_{1},  \dots, x_{n}$ be elements  of $H$ such  that their images
$v_{1}, \dots , v_{n}$ are a basis  of $V = H/H'$. Now if in $\Span{u}
\oplus V$ we have
\begin{equation*}
  u \phi = u + \sigma_{1} v_{1} + \dots + \sigma_{n} v_{n},
\end{equation*}
where $\sigma = [\sigma_{1}, \dots, \sigma_{n}]$ with respect to
the basis of the $v_{i}$, then
\begin{align*}
  [u, v_{1}] 
  &= [u, v_{1}] \phi = [u \phi, v_{1}]
  \\&=
  [u + \sigma_{1} v_{1} + \dots + \sigma_{n} v_{n}, v_{1}]
  \\&=
  [u, v_{1}] + \sigma_{2} [v_{2}, v_{1}] + \dots + \sigma_{n} [v_{n}, v_{1}],
\end{align*}
so that
\begin{equation*}
  \sigma_{2} [v_{2}, v_{1}] + \dots + \sigma_{n} [v_{n}, v_{1}] = 0.
\end{equation*}
Now the elements $ v_{2} \wedge v_{1}, \dots , v_{n} \wedge v_{1}$ are
independent in $\Lambda^{2} V$. Since  $K = \Set{0}$, we have $G' \cong
\Lambda^{2}  V$,  so that  $\sigma_{2}  =  \dots  = \sigma_{n}  =  0$.
Considering  $[u,  v_{2}]  \phi  =  [u,  v_{2}]$,  we  see  that  also
$\sigma_{1} = 0$, so that  $\sigma = 0$, the matrix in~\eqref{eq:mat2}
is the identity,  and thus $\phi$ induces the  identity modulo $Z(G) =
\Phi(G)$.

Using   once  more   the   approach  of
Section~\ref{sec:Zurek} (or appealing again to~\cite[Theorem
  4]{AdYen}),  we  have  that 
$\Aut(G) =  \Aut_{c}(G)$  is 
isomorphic to $(\Hom(G, Z(G))$ under the circle operation. Take two
elements $\gamma_{1}, 
\gamma_{2} \in \Hom(G, Z(G))$, which are easily seen to be of the form
\begin{equation*}
  \gamma_{k} : y \mapsto y^{p s_{k}} \pmod{G'}, 
  \qquad
  x_{i} \mapsto 1 \pmod{G'}, \text{ for
      all $i$},
\end{equation*}
for some $s_{1}, s_{2}$. Since $G' \le
\ker(\gamma_{k})$, and $(G')^{p} = 1$, we have
\begin{equation*}
  y \gamma_{1} \gamma_{2} = y^{p^{2} s_{1} s_{2}},
  \qquad
  x_{i} \gamma_{1} \gamma_{2}  = 1, \text{ for all
      $i$},
\end{equation*}
from which it follows that $\Aut(G) = \Aut_{c}(G)$ is abelian. 
Moreover $G' < \Span{G', y^{p}} = \Phi(G) = Z(G)$.

\providecommand{\bysame}{\leavevmode\hbox to3em{\hrulefill}\thinspace}
\providecommand{\MR}{\relax\ifhmode\unskip\space\fi MR }
\providecommand{\MRhref}[2]{%
  \href{http://www.ams.org/mathscinet-getitem?mr=#1}{#2}
}
\providecommand{\href}[2]{#2}

\end{document}